\documentclass[11pt,twoside]{article}
\usepackage{multicol}
\usepackage[greek,english]{babel}
\usepackage[vcentering,dvips]{geometry}
\usepackage{makeidx}
\usepackage{amsfonts}
\usepackage{amssymb}
\usepackage{eucal}
\usepackage{amsxtra}
\usepackage{color}
\usepackage{setspace}
\usepackage{times}
\usepackage{titlesec}
\usepackage{sectsty}
\usepackage{fancyhdr}
\usepackage{fancyheadings}
\allsectionsfont{\large}
\usepackage{url}
\usepackage{doi}
\usepackage{hyperref}
\hypersetup{
    colorlinks,%
    citecolor=black,%
    filecolor=black,%
    linkcolor=black,%
    urlcolor=black
}
\usepackage[T1]{fontenc}
\usepackage{selinput}
\SelectInputMappings{%
  adieresis={ä},
  eacute={é},
  Lcaron={Ľ},
}

\addtocounter{section}{0} \numberwithin{equation}{section}
\newtheorem{theorem}{Theorem}[section]
\newtheorem{proposition}[theorem]{Proposition}
\newtheorem{lemma}[theorem]{Lemma}
\newtheorem{remark}[theorem]{Remark}
\newtheorem{corollary}[theorem]{Corollary}

\newcommand{\MN}[1][]{\ensuremath{{M_{n}(c)}\;}}
\newcommand{\MNc}[1][]{\ensuremath{{M_{2}(c)}\;}}

\newcommand{\M}{{\emph{M }}}
\newcommand{\K}{{\emph{P }}}

\hypersetup{urlcolor=blue }

\begin{document}
\title{\textbf{Nullity distributions on real hypersurfaces in non-flat complex space forms}}
\author{\textsc{Konstantina Panagiotidou}}
\date{}
\maketitle

\begin{abstract}
In this paper the result of real hypersurfaces in non-flat complex space forms, whose structure vector field $\xi$ belongs to the $\kappa$-nullity distribution is extended in case of three dimensional real hypersurfaces in non-flat complex space forms. Furthermore, generalization of notion ($\kappa$,$\mu$)-nullity distribution defined on real hypersurfaces and results of real hypersurfaces, whose structure vector field $\xi$ belongs to the previous distribution are provided. Finally, the notion of ($\kappa$,$\mu$,$\nu$)-nullity distribution is introduced in case of real hypersurfaces in non-flat complex space forms and real hypersurfaces, whose structure vector field $\xi$ belongs to the previous distribution are studied.
\end{abstract}

\noindent\small{\emph{Keywords}: Real hypersurfaces, Non-flat Complex Space Forms, Nullity Distributions, Structure Vector Field.\\}

\rhead[\centering{K. Panagiotidou}]{\thepage}
\lhead[\thepage]{\centering{Nullity distributions}}

\section{\textsc{Introduction}}
A \emph{complex space form} is an $n$-dimensional Kähler manifold of constant holomorphic sectional curvature \emph{c}. A complete and simply connected complex space form is complex analytically isometric to complex projective space $\mathbb{C}P^{n}$ if $c>0$, or to complex Euclidean space $\mathbb{C}^{n}$ if $c=0$, or to complex hyperbolic space $\mathbb{C}H^{n}$ if $c<0$. The complex projective and complex hyperbolic spaces are called \emph{non-flat complex space forms}, since $c\neq0$ and the symbol $M_{n}(c)$ is used to denote them when it is not necessary to distinguish them.

A real hypersurface \M is an immersed submanifold with real co-dimension one in $M_{n}(c)$. The Kähler structure ($J, G$), where $J$ is the complex structure and $G$ is the Kähler metric of $M_{n}(c)$, induces on \M an almost contact metric structure ($\varphi,\xi,\eta,g$). The vector field $\xi$ is called \emph{structure vector field} and when it is an eigenvector of the shape operator $A$ of \M the real hypersurface is called \emph{Hopf hypersurface} with corresponding eigenvalue is $\alpha=g(A\xi,\xi)$.

The study of real hypersurfaces \M in $M_{n}(c)$ was initiated by Takagi, who classified homogeneous real hypersurfaces in $\mathbb{C}P^{n}$ and divided them into six types, namely ($A_{1}$), ($A_{2}$), ($B$), ($C$), ($D$) and ($E$) (\cite{T2}, \cite{T1}). These real hypersurfaces are Hopf ones with constant principal curvatures. In case of $\mathbb{C}H^{n}$, the study of real hypersurfaces with constant principal curvatures  was started by Montiel in \cite{Mo} and completed by Berndt in \cite{Ber}. They are divided into two types, namely ($A$) and ($B$), depending on the number of constant principal curvatures and they are homogeneous and Hopf hypersurfaces.

Many geometers have studied real hypersurfaces in non-flat complex space forms when certain geometric conditions are satisfied. An important condition is that of the shape operator $A$ commuting with the structure tensor $\varphi$. More precisely, the following Theorem owed to Okumura in case of $\mathbb{C}P^{n}$  (\cite{Ok}) and to Montiel and Romero in case of $\mathbb{C}H^{n}$ (\cite{MR}) plays an important role in the proof of other Theorems.

\begin{theorem}\label{Th-Ok-Mo-Ro}
Let M be a real hypersurface of $M_{n}(c)$, $n\geq2$. Then $A\varphi=\varphi A$, if and only if M is locally congruent to a homogeneous real hypersurface of type (A). More precisely\\
 in case of $\mathbb{C}P^{n}$\\
    $(A_{1})$   a geodesic hypersphere of radius r , where $0<r<\frac{\pi}{2}$,\\
    $(A_{2})$  a tube of radius r over a totally geodesic $\mathbb{C}P^{k}$,$(1\leq k\leq n-2)$, where $0<r<\frac{\pi}{2}.$\\
 In case of $\mathbb{C}H^{n}$\\
    $(A_{0})$   a horosphere in $ \mathbb{C}H^{n}$, i.e a Montiel tube,\\
    $(A_{1})$  a geodesic hypersphere or a tube over a totally geodesic complex hyperbolic hyperplane $\mathbb{C}H^{n - 1}$,\\
    $(A_{2}) $  a tube over a totally geodesic $\mathbb{C}H^{k}$ $(1\leq k\leq n-2)$.
\end{theorem}
In \cite{TA} Tanno introduced the notion of $\kappa$ - \emph{nullity distribution} for Riemannian manifolds,
\[N(\kappa):P\rightarrow N_{P}(\kappa),\]
and $N_{P}(\kappa)$ is given by
\[N_{P}(\kappa)=\{Z\;\in\;T_{P}M:R(X,Y)Z=\kappa[g(Y,Z)X-g(X,Z)Y]\},\;\;\mbox{for any $X$, $Y$ $\in$ $T_{P}M$}.\]

In case of real hypersurfaces in non-flat complex space forms, in \cite{CH1} and \cite{CH2} Cho studied real hypersurfaces in \MN, $n\geq3$, whose structure vector field $\xi$ belongs to $\kappa$-nullity distribution with $\kappa$ smooth function. It was proved that such real hypersurfaces are of type ($A$) and $\kappa$ is constant.

In \cite{CK} the notion of ($\kappa$, $\mu$) - \emph{nullity distribution} was introduced in the following way
\[N(\kappa, \mu):P\rightarrow N_{P}(\kappa, \mu),\;\;\mbox{with $(\kappa, \mu)$ $\in$ $\mathbb{R}^{2}$} \]
and $N_{P}(\kappa, \mu)$ is given by
\[N_{P}(\kappa, \mu)=\{Z\;\in\;T_{P}M:R(X,Y)Z=(\kappa I+\mu A)[g(Y,Z)X-g(X,Z)Y]\},\mbox{for any $X$, $Y$ $\in$ $T_{P}M$}.\]
In the above relation $I$ denotes the identity and $A$ the shape operator of real hypersurface. Furthermore, in \cite{CK} it was proved that real hypersurfaces in non-flat complex space forms, whose structure vector field $\xi$ belongs to ($\kappa$,$\mu$)-nullity distribution with ($\kappa,\mu$) $\in$ $\mathbb{R}^{2}$ are Hopf and classification in case of $(0,\mu)$-nullity distribution and $(\kappa,0)$-nullity distribution with $\alpha\neq0$ is obtained.

In \cite{KMP} the notion of ($\kappa$,$\mu$,$\nu$)-nullity distribution was introduced and studied for contact metric manifolds. Motivated by their work, in this paper the notion of ($\kappa$, $\mu$, $\nu$) - \emph{nullity distribution} is introduced for real hypersurfaces in \MN, $n\geq2$, in the following way
\[N(\kappa, \mu, \nu):P\rightarrow N_{P}(\kappa, \mu, \nu),\;\;\mbox{with $\kappa$, $\mu$, $\nu$ smooth functions}\]
 and  $N_{P}(\kappa, \mu, \nu)$ is given by
\[N_{P}(\kappa, \mu ,\nu)=\{Z\;\in\;T_{P}M:R(X,Y)Z=\kappa[\eta(Y)X-\eta(X)Y]+\mu[\eta(Y)AX-\eta(X)AY]\]
\[+\nu[\eta(Y)\varphi AX-\eta(X)\varphi AY]\},\;\;\mbox{for any $X$, $Y$ $\in$ $T_{P}M$}.\]

Motivated by the work that so far has been done, the following questions raised naturally

\textbf{Questions}: 1) \emph{Do there exist real hypersurfaces in non-flat complex space forms, whose structure vector field $\xi$ belongs to ($\kappa$, $\mu$) - nullity distribution, with ($\kappa$, $\mu$) smooth functions?}\\
2) \emph{Do there exist real hypersurfaces in non-flat complex space forms, whose structure vector field $\xi$ belongs to ($\kappa$, $\mu$, $\nu$) - nullity distribution, with $\kappa$, $\mu$, $\nu$ smooth functions?}

The aim of this paper is first to extend the results in \cite{CH1} and \cite{CH2} in case of three dimensional real hypersurfaces in \MNc. More precisely, the following Theorem is proved.
\begin{theorem}\label{tk}
Every real hypersurface M in \MNc, whose structure vector field $\xi$ belongs to $\kappa$-nullity distribution is locally congruent either to a real hypersurface of type ($A$) with $\kappa$ constant or to a real hypersurface with $A\xi=0$ and $\kappa$ constant.
\end{theorem}
Next, the first question mentioned above is answered in case of $\kappa$, $\mu$ are non-constant smooth functions
\begin{theorem}\label{tkn}
There do not exist real hypersurfaces M in \MN, $n\geq2$, whose structure vector field $\xi$ belongs to ($\kappa$, $\mu$) - nullity distribution with $\kappa$,$\mu$ non-constant smooth functions.
\end{theorem}
Finally, the following Theorem provides an answer in the second question in case of $\kappa$,$\mu$,$\nu$ are non-constant smooth functions.
\begin{theorem}\label{tknm}
There do not exist real hypersurfaces M in \MN, $n\geq2$, whose structure vector field $\xi$ belongs to ($\kappa$, $\mu$, $\nu$) - nullity distribution with $\kappa$, $\mu$, $\nu$ non-constant smooth functions.
\end{theorem}
This paper is organized as follows: in Section 2 relations and basic results which hold for real hypersurfaces in non-flat complex space forms are presented. In Section 3 the proof of Theorem \ref{tk}. In Section 4 the proof of Theorem \ref{tkn}, which generalizes the results obtained in \cite{KMP} is included. Finally, in Section 5 Theorem \ref{tknm} is proved and at the end of Section open problems for further research are provided.

\section{\textsc{Preliminaries}}
Throughout this paper all manifolds, vector fields etc are assumed to be of class $C^{\infty}$ and all manifolds are assumed to be connected. Furthermore, the real hypersurfaces $M$ are supposed to be without boundary.

Let \M be a real hypersurface immersed in $(M_{n}(c),G)$ with complex structure $J$ of constant holomorphic sectional curvature $c$. Let $N$ be a unit normal vector field on \M and $\xi=-JN$ the structure vector field of \M. For any vector field $X$ tangent to \M relation
\[JX=\varphi X+\eta(X)N\]
holds, where $\varphi X$ and $\eta(X)N$ are respectively the tangential and the normal
component of $JX$. The Riemannian connections
$\overline{\nabla}$ in $M_{n}(c)$ and $\nabla$ in $M$ are related
for any vector fields $X$, $Y$ on $M$ by
\[\overline{\nabla}_{X}Y=\nabla_{X}Y+g(AX,Y)N,\]
where $g$ is the Riemannian metric induced from the metric $G$.

The shape operator $A$ of the real hypersurface \M in $M_{n}(c)$ with respect to $N$ is given by
\[\overline{\nabla}_{X}N=-AX.\]
The real hypersurface \M has an almost contact metric structure $(\varphi,\xi,\eta, g)$ induced from $J$ on $M_{n}(c)$, where $\varphi$  is the \emph{structure tensor} which is a tensor field of type (1,1) and $\eta$ is an
1-form on \M such that
\[g(\varphi X,Y)=G(JX,Y),\hspace{20pt}\eta(X)=g(X,\xi)=G(JX,N).\]
Moreover, the following relations hold
\begin{eqnarray}\label{eq-1}
\varphi^{2}X=-X+\eta(X)\xi,\hspace{20pt}
\eta\circ\varphi=0,\hspace{20pt} \varphi\xi=0,\hspace{20pt}
\eta(\xi)=1,\nonumber\\
\hspace{20pt}
g(\varphi X,\varphi
Y)=g(X,Y)-\eta(X)\eta(Y),\hspace{10pt}g(X,\varphi Y)=-g(\varphi
X,Y).\label{eq-2}\nonumber\
\end{eqnarray}
The fact that $J$ is parallel implies $\bar{\nabla} J=0$. The last relation leads to
\begin{eqnarray}\label{eq-3}
\nabla_{X}\xi=\varphi
AX,\hspace{20pt}(\nabla_{X}\varphi)Y=\eta(Y)AX-g(AX,Y)\xi.
\end{eqnarray}
    The ambient space $M_{n}(c)$ is of constant holomorphic sectional
curvature $c$. Thus, the Gauss and Codazzi equations to are respectively given by
\begin{eqnarray}\label{eq-4}
R(X,Y)Z=\frac{c}{4}[g(Y,Z)X-g(X,Z)Y+g(\varphi Y ,Z)\varphi
X
\end{eqnarray}
$$-g(\varphi X,Z)\varphi Y-2g(\varphi X,Y)\varphi
Z]+g(AY,Z)AX-g(AX,Z)AY,$$
\begin{eqnarray}\label{eq-5}
\hspace{10pt}
(\nabla_{X}A)Y-(\nabla_{Y}A)X=\frac{c}{4}[\eta(X)\varphi
Y-\eta(Y)\varphi X-2g(\varphi X,Y)\xi],
\end{eqnarray}
where $R$ is the Riemannian curvature tensor on \M and $X$, $Y$, $Z$ are any vector fields on \M.

At every point $P$ $\in$ \M the tangent space $T_{P}M$ is decomposed
\begin{eqnarray}
T_{P}M=span\{\xi\}\oplus \mathbb{D},\nonumber\
\end{eqnarray}
where $\mathbb{D}=\ker\eta=\{X\;\;\in\;\; T_{P}M:\eta(X)=0\}$ and is called (\emph{maximal}) \emph{holomorphic distribution}, (if $n\geq3$).
The above relation implies that the vector field $A\xi$ can be written
 \begin{eqnarray}\label{eq-7}
 A\xi=\alpha\xi+\beta U,\nonumber\
 \end{eqnarray}
 where $\beta=|\varphi\nabla_{\xi}\xi|$ and
 $U=-\frac{1}{\beta}\varphi\nabla_{\xi}\xi\;\in\;\ker(\eta)$ is a unit vector field, provided
 that $\beta\neq0$.

The following Theorem is owed to Maeda in  case of $\mathbb{C}P^{n}$ \cite{Maeda} and to Montiel \cite{Mo}  in case of $\mathbb{C}H^{n}$(also Corollary 2.3 in \cite{NR1}).
\begin{theorem}\label{Ma-Mo}
Let M be a Hopf hypersurface in $M_{n}(c)$, $n\geq2$. Then\\
i) $\alpha$ is constant.\\
ii) If $W$ is a vector field which belongs to $\mathbb{D}$ such that $AW=\lambda_{1} W$, then
\begin{eqnarray}\label{eq-A}
(\lambda_{1}-\frac{\alpha}{2})A\varphi W=(\frac{\lambda_{1}\alpha}{2}+\frac{c}{4})\varphi W.\nonumber\
\end{eqnarray}
iii) If the vector field $W$ satisfies $AW=\lambda_{1} W$ and $A\varphi W=\lambda_{2}\varphi W$ then
\begin{eqnarray}\label{eq-B}
\lambda_{1}\lambda_{2}=\frac{\alpha}{2}(\lambda_{1}+\lambda_{2})+\frac{c}{4}.
\end{eqnarray}
\end{theorem}

\begin{remark}\label{di}
In case of real hypersurfaces of dimension greater than three the third case of Theorem \ref{Ma-Mo} occurs when $\alpha^{2}+c\neq0$, since in this case relation $\lambda_{1}\neq\frac{\alpha}{2}$ holds.
\end{remark}

\begin{remark}\label{three}
In case of three dimensional Hopf hypersurfaces it can be always considered a local orthonormal basis $\{W,\varphi W, \xi\}$ at some point $P$ $\in$ M such that $AW=\lambda_{1} W$ and $A\varphi W=\lambda_{2}\varphi W$.
\end{remark}

\subsection{Auxiliary facts about three dimensional real hypersurfaces in complex space forms}
Let \M be a real hypersurface in $M_{2}(c)$ with local orthonormal basis $\{U, \varphi U, \xi\}$ at some point \K of $M$. \begin{lemma}\label{lemma-1}
Let M be a non-Hopf real hypersurface in $M_{2}(c)$. The following relations hold on M
\begin{eqnarray}
&&AU=\gamma U+\delta\varphi U+\beta\xi,\hspace{20pt} A\varphi U=\delta U+\rho\varphi U,\hspace{20pt}A\xi=\alpha\xi+\beta U,\label{eq-8}\\
&&\nabla_{U}\xi=-\delta U+\gamma\varphi U,\hspace{20pt}
\nabla_{\varphi U}\xi=-\rho U+\delta\varphi U,\hspace{20pt}
\nabla_{\xi}\xi=\beta\varphi U,\label{eq-9}\nonumber\\
&&\nabla_{U}U=\kappa_{1}\varphi U+\delta\xi,\hspace{20pt}
\nabla_{\varphi U}U=\kappa_{2}\varphi U+\rho\xi,\hspace{20pt}
\nabla_{\xi}U=\kappa_{3}\varphi U,\label{eq-10}\nonumber\\
&&\nabla_{U}\varphi U=-\kappa_{1}U-\gamma\xi,\hspace{5pt}
\nabla_{\varphi U}\varphi U=-\kappa_{2}U-\delta\xi,\hspace{5pt}
\nabla_{\xi}\varphi U=-\kappa_{3}U-\beta\xi,\label{eq-11}\nonumber\
\end{eqnarray}
where $\alpha$, $\beta$, $\gamma,\delta,\rho,\kappa_{1},\kappa_{2},\kappa_{3}$ are
smooth functions on M and $\beta\neq0$.
\end{lemma}

\begin{remark}
The proof of Lemma \ref{lemma-1} is included in \cite{PX3}.
\end{remark}

The Codazzi equation (\ref{eq-5}) for $X$ $\in$ $\{U, \varphi U\}$ and $Y=\xi$ because of Lemma \ref{lemma-1} implies the following relations
\begin{eqnarray}
\xi\delta&=&\alpha\gamma+\beta\kappa_{1}+\delta^{2}+\rho\kappa_{3}+\frac{c}{4}-\gamma\rho-\gamma\kappa_{3}-\beta^{2},\label{6}\\
(\varphi U)\alpha&=&\alpha\beta+\beta\kappa_{3}-3\beta\rho,\label{9}\\
(\varphi U)\beta&=&\alpha\gamma+\beta\kappa_{1}+2\delta^{2}+\frac{c}{2}-2\gamma\rho+\alpha\rho,\label{10}
\end{eqnarray}
and for $X=U$ and $Y=\varphi U$
\begin{eqnarray}
U\delta-(\varphi U)\gamma&=&\rho\kappa_{1}-\kappa_{1}\gamma-\beta\gamma-2\delta\kappa_{2}-2\beta\rho.\label{11}
\end{eqnarray}

\section{$\kappa$-Nullity Distribution}
Let \M be a real hypersurface in \MNc, whose structure vector fields $\xi$ belongs to $\kappa$ - nullity distribution, i.e.
\begin{eqnarray}\label{bs-1}
R(X,Y)\xi=\kappa[\eta(Y)X-\eta(X)Y],\;\;\; \mbox{where $\kappa$ is a smooth function and $X$,$Y$ $\in$ $TM$.}
\end{eqnarray}
Let $\mathcal{N}$  be the open subset of $M$ such that
\[\mathcal{N}=\{P\;\;\in\;\;M:\;\beta\neq0\;\;\mbox{in a neighborhood of $P$}\}.\]
Relation (\ref{bs-1}) for $X=U$ and $Y=\xi$ because of (\ref{eq-4}) and (\ref{eq-8}) implies
\begin{eqnarray}\label{k1}
\delta=0\;\;\mbox{and}\;\;\kappa=\frac{c}{4}+\alpha\gamma-\beta^{2}
\end{eqnarray}
and for $X=U$ and $Y=\varphi U$ due to (\ref{eq-4}), (\ref{eq-8}) and $\delta=0$ implies $\rho=0$. Furthermore, relation (\ref{bs-1}) for $X=\varphi U$ and $Y=\xi$ due to (\ref{eq-4}), (\ref{eq-8}) and $\delta=\rho=0$ yields $\kappa=\frac{c}{4}$. Combination of the second of (\ref{k1}) with the last one results in $\alpha\gamma=\beta^{2}$. Differentiation of the last one with respect to $\varphi U$ taking into account relations (\ref{6}), (\ref{9}), (\ref{10}), (\ref{11}), $\alpha\gamma=\beta^{2}$ and $\delta=\rho=0$ implies $c=0$, which is a contradiction.

Therefore, $\mathcal{N}$ is empty and the following proposition has been proved
\begin{proposition}
Every real hypersurface M in \MNc, whose structure vector field $\xi$ satisfies (\ref{bs-1}) is Hopf.
\end{proposition}

Since \M is a Hopf hypersurface Theorem \ref{Ma-Mo} and remark \ref{three} hold. Thus, relation (\ref{bs-1}) for $X=W$ and $Y=\xi$ and for $X=\varphi W$ and $Y=\xi$, owing to (\ref{eq-4}), $AW=\lambda_{1} W$ and $A\varphi W=\lambda_{2}\varphi W$ yields respectively
\[\kappa=\frac{c}{4}+\alpha\lambda_{1}\;\;\mbox{and}\;\;\kappa=\frac{c}{4}+\alpha\lambda_{2}.\]
Combination of the latter results in
\[\alpha(\lambda_{1}-\lambda_{2})=0.\]

If $\alpha=0$ then in case of $\mathbb{C}P^{2}$, \M is locally congruent to a non-homogeneous real hypersurface considered as a tube of radius $r=\frac{\pi}{4}$ over a holomorphic curve or to a geodesic hypersphere of radius $r=\frac{\pi}{4}$. In case of $\mathbb{C}H^{2}$, \M is locally congruent to a Hopf hypersurface with $A\xi=0$ (for the construction of such real hypersurfaces see \cite{IR1}).

If $\alpha\neq0$ then $\lambda_{1}=\lambda_{2}$  which implies
\[(A\phi-\phi A)X=0,\;\;\mbox{for any $X$ tangent to \M}.\]
So due to Theorem {\ref{Th-Ok-Mo-Ro} \M is locally congruent to a real hypersurface of type ($A$).

Conversely, if \M is a real hypersurface of type ($A$), then the shape operator is given by
\[A\xi=\alpha\xi\;\;\mbox{and}\;\;AW=\lambda W,\;\;\mbox{for any $W$ $\in$ $\mathbb{D}$ and $\alpha$, $\lambda$ constants}.\]

Combination of (\ref{eq-4}) with (\ref{bs-1}) for any $X=W$ $\in$ $\mathbb{D}$ and $Y=\xi$ due to $AW=\lambda W$ implies that the structure vector field $\xi$ belongs to $\kappa$-nullity distribution when
\[\kappa=\frac{c}{4}+\alpha\lambda.\]
Thus, in case of $\mathbb{C}P^{2}$ we have $c=4$ and when \M is locally congruent to geodesic hypersphere then $\kappa=\cot^{2}(r)$.

In case of $\mathbb{C}H^{2}$ we have $c=-4$ and
\begin{itemize}
\item when \M is locally congruent to a horosphere, then $\kappa=1$
\item when \M is locally congruent to a geodesic hypersphere, then $\kappa=\coth^{2}(r)$,
\item when \M is locally congruent to a tube over totally geodesic $\mathbb{C}H^{1}$, then $\kappa=\tanh^{2}(r)$.
\end{itemize}

If \M is a Hopf hypersurface with $A\xi=0$ then following similar steps as above the structure vector field $\xi$ belongs to $\kappa$-nullity distribution when
\[\kappa=\frac{c}{4}.\]
So, in case of $\mathbb{C}P^{2}$ we have $\kappa=1$ and in case of $\mathbb{C}H^{2}$ we have $\kappa=-1$.

\section{($\kappa$, $\mu$)-Nullity Distribution}

Let $M$ be a real hypersurface in $M_{n}(c)$, $n\geq2$, whose structure vector field $\xi$ belongs to $(\kappa, \mu)$- nullity distribution, i.e.
\begin{eqnarray}\label{bs-2}
R(X,Y)\xi=\kappa[\eta(Y)X-\eta(X)Y]+\mu[\eta(Y)AX-\eta(X)AY],
\end{eqnarray}
with $\kappa$, $\mu$ non-constant smooth functions.

Consider $\mathcal{N}$ the open subset of \M such that
\[\mathcal{N}=\{P\;\;\in\;\;M:\;\beta\neq0\;\;\mbox{in a neighborhood of $P$}\}.\]

On $\mathcal{N}$ relation (\ref{bs-2}) for $X=U$ and $Y=\varphi U$ because of (\ref{eq-4}), $A\xi=\alpha\xi+\beta U$ and $\beta\neq0$ implies $A\varphi U=0$. Furthermore, relation (\ref{bs-2}) for $X=\varphi U$ and $Y=\xi$ due to the last relation yields $\kappa=\frac{c}{4}$.

Relation $A\varphi U=0$ results in $g(AU, \varphi U)=g(A\varphi U, U)=0$. So $AU$ can be written as $AU=\gamma U+\beta\xi+t Z$, where $Z$ is a unit vector field in $\mathbb{D}_{U}=span\{U, \varphi U, \xi\}^{\perp}$. Moreover, relation (\ref{bs-2}) for $X=U$ and $Y=\xi$, owing to (\ref{eq-4}), the last one and $\kappa=\frac{c}{4}$ yields
\[(\mu\gamma+\beta^{2}-\alpha\gamma)U+\mu\beta\xi+(\mu t-\alpha t)Z=0.\]
The inner product of the above relation with $\xi$ due to $\beta\neq0$ implies $\mu=0$ and relation (\ref{bs-2}) becomes
\[R(X,Y)\xi=\kappa[\eta(Y)X-\eta(X)Y].\]
Thus, $\xi$ belongs to a $\kappa$-nullity distribution and because of Main Theorem in \cite{CH1}, Lemma 2 in \cite{CH2} and Theorem \ref{tk} of the present paper it is proved that $\mathcal{N}$ is empty. Thus,

\begin{proposition}\label{proposition-2*}
Every real hypersurface in $M_{n}(c)$, $n\geq2$ whose structure vector field $\xi$ satisfies relation (\ref{bs-2}) is Hopf.
\end{proposition}

Since \M is a Hopf hypersurface in $M_{n}(c)$, $n\geq2$, we consider two cases

\textbf{Case I:} $\alpha^{2}+c\neq0$.

Consider a vector field $W$ $\in$ $\mathbb{D}$ such that $AW=\lambda_{1}W$. Then Theorem \ref{Ma-Mo} and remark \ref{di} hold. So, relation $A\varphi W=\lambda_{2}\varphi W$ holds.

So, relation (\ref{bs-2}) for $X=W$ and $Y=\xi$ and for $X=\varphi W$ and $Y=\xi$ because of relation (\ref{eq-4}) respectively yields
\[\kappa+\mu\lambda_{1}=\frac{c}{4}+\alpha\lambda_{1}\;\;\mbox{and}\;\;\kappa+\mu\lambda_{2}=\frac{c}{4}+\alpha\lambda_{2}.\]
Combination of the above relations implies
\[(\lambda_{1}-\lambda_{2})(\alpha-\mu)=0.\]

Suppose that $\lambda_{1}\neq\lambda_{2}$ then $\mu=\alpha$. Substitution of the latter in $\kappa+\mu\lambda_{1}=\frac{c}{4}+\alpha\lambda_{1}$ results in $\kappa=\frac{c}{4}$, which is a contradiction.

Therefore, on \M $\lambda_{1}=\lambda_{2}$ which implies
\[(A\varphi-\varphi A)X=0,\;\;\mbox{for any $X$ tangent to \M}.\]
So, because of Theorem \ref{Th-Ok-Mo-Ro} \M is locally congruent to a real hypersurface of type ($A$).

\textbf{Case II: }$\alpha^{2}+c=0$.

In this case the ambient space is $\mathbb{C}H^{n}$, $n\geq2$ and the above relation implies that $\alpha\neq0$. First suppose that $\lambda_{1}\neq\frac{\alpha}{2}$. Then relation (\ref{eq-B}) yields $\lambda_{2}=\frac{\alpha}{2}$. Following similar steps as in previous case we obtain
\[(\lambda_{1}-\frac{\alpha}{2})(\alpha-\mu)=0.\]
Since, $\lambda_{1}\neq\frac{\alpha}{2}$, then $\mu=\alpha$, which is contradiction.

So $\lambda_{1}=\frac{\alpha}{2}$ is the only eigenvalue in $\mathbb{D}$ and \M is locally congruent to a horosphere.

Conversely, if \M is a real hypersurface of type ($A$) then \M has either two or three constant principal curvatures and the maximal holomorphic distribution is $\varphi$-invariant. First, suppose that \M has two constant principal curvatures, then the shape operator is given by
\[A\xi=\alpha\xi\;\;\mbox{and}\;\;AW=\lambda W,\;\;\mbox{for any W $\in$ $\mathbb{D}$}.\]

Then the structure vector field $\xi$ belongs to ($\kappa$,$\mu$)-nullity distribution when
\[\kappa=\frac{c}{4}\;\;\mbox{and}\;\;\mu=\alpha.\]
Indeed combination of relation (\ref{eq-4}) with (\ref{bs-2}) because of the form of the shape operator yields
\[\frac{c}{4}+\alpha\lambda=\kappa+\mu\lambda.\]
The above two polynomials of $\lambda$ are equal when $\kappa=\frac{c}{4}$ and $\mu=\alpha$, which is a contradiction.

Finally, suppose that \M has three distinct constant principal curvatures. Then the shape operator is given by
\[A\xi=\alpha\xi,\;\;AW_{1}=t_{1}W\;\;\mbox{and}\;\;AW_{2}=t_{2}W_{2},\;\;W_{1},W_{2}\;\in\;\mathbb{D}.\]
Then the structure vector field $\xi$ belongs to ($\kappa$,$\mu$)-nullity distribution when
\[\kappa=\frac{c}{4}\;\;\mbox{and}\;\;\mu=\alpha.\]
Indeed combination of relation (\ref{eq-4}) with (\ref{bs-2}) due to the form of the shape operator implies
\[\frac{c}{4}+\alpha t_{1}=\kappa+\mu t_{1}\;\;\mbox{and}\;\;\frac{c}{4}+\alpha t_{2}=\kappa+\mu t_{2}.\]
Combining the last two relation and taking into account that $t_{1}\neq t_{2}$ results in $\kappa=\frac{c}{4}\;\;\mbox{and}\;\;\mu=\alpha$, which is a contradiction and this completes the proof of Theorem \ref{tkn}.

\section{($\kappa$, $\mu$, $\nu$) - Nullity Distribution}
Let $M$ be a real hypersurface in $M_{n}(c)$, $n\geq2$, whose structure vector field $\xi$ belongs to $(\kappa, \mu, \nu)$- nullity distribution, i.e.
\begin{eqnarray}\label{bs-3}
R(X,Y)\xi=\kappa[\eta(Y)X-\eta(X)Y]+\mu[\eta(Y)AX-\eta(X)AY]+\nu[\eta(Y)\varphi AX-\eta(X)\varphi AY],
\end{eqnarray}
where $\kappa$, $\mu$, $\nu$ are non-constant smooth functions.

Let $\mathcal{N}$ be the open subset of \M such that
\[\mathcal{N}=\{P\;\;\in\;\;M:\;\beta\neq0\;\;\mbox{in a neighborhood of $P$}\}.\]

On $\mathcal{N}$ relation (\ref{bs-3}) for $X=U$ and $Y=\varphi U$ because of (\ref{eq-4}), $A\xi=\alpha\xi+\beta U$ and $\beta\neq0$ implies $A\varphi U=0$. Furthermore, relation (\ref{bs-3}) for $X=\varphi U$ and $Y=\xi$ due to the last relation yields $\kappa=\frac{c}{4}$.

Relation $A\varphi U=0$ results in $g(AU, \varphi U)=g(A\varphi U, U)=0$. So $AU$ can be written as $AU=\gamma U+\beta\xi+t Z$, where $Z$ is a unit vector field in $\mathbb{D}_{U}=span\{U, \varphi U, \xi\}^{\perp}$. Moreover, relation (\ref{bs-3}) for $X=U$ and $Y=\xi$, owing to (\ref{eq-4}), the last one and $\kappa=\frac{c}{4}$ yields
\begin{eqnarray}\label{t1}
(\mu\gamma+\beta^{2}-\alpha\gamma)U+\nu\gamma\varphi U+\mu\beta\xi+(\mu t-\alpha t)Z+\nu t\varphi Z=0.
\end{eqnarray}
The inner product of the above relation with $\xi$ since $\beta\neq0$ implies $\mu=0$ and with $\varphi U$ yields $\nu\gamma=0$. If $\nu\neq0$ then $\gamma=0$. The inner product of (\ref{t1}) with $U$ because of $\mu=\gamma=0$ results in $\beta=0$, which is a contradiction.

So on $\mathcal{N}$ relation $\nu=0$ holds and relation (\ref{bs-3}) becomes
\[R(X,Y)\xi=\kappa[\eta(Y)X-\eta(X)Y].\]
Thus, $\xi$ belongs to a $\kappa$-nullity distribution and because of Main Theorem in \cite{CH1},Lemma 2 in \cite{CH2} and Theorem \ref{tk} of the present paper it is concluded that $\mathcal{N}$ is empty and the following Proposition has been proved

\begin{proposition}\label{proposition-1*}
Every real hypersurface in $M_{n}(c)$, $n\geq2$, whose structure vector field $\xi$ satisfies relation (\ref{bs-3}) is Hopf.
\end{proposition}

Since \M is a Hopf hypersurface in $M_{n}(c)$, $n\geq2$, two cases are considered

\textbf{Case I:} $\alpha^{2}+c\neq0$.

Let $W$ be a vector field which belongs to $\mathbb{D}$ such that $AW=\lambda_{1}W$. In this case $\lambda_{1}\neq\frac{\alpha}{2}$, so $A\varphi W=\lambda_{2}\varphi W$ and relations of Theorem \ref{Ma-Mo} and remark \ref{di} hold.

The inner product of relation (\ref{bs-3}) for $X=W$ and $Y=\xi$ with $\varphi W$ and for $X=\varphi W$ and $Y=\xi$ with $W$ because of (\ref{eq-4}) and the above relations respectively yields
\[\nu\lambda_{1}=0\;\;\mbox{and}\;\;\nu\lambda_{2}=0.\]
Combination of the last two relations results in
\[\nu(\lambda_{1}-\lambda_{2})=0.\]

Suppose that $\nu\neq0$ then $\lambda_{1}=\lambda_{2}$ and relation $\nu\lambda_{1}=0$ results in $\lambda_{1}=\lambda_{2}=0$. Substitution of the latter in (\ref{eq-B}) implies $c=0$, which is a contradiction.

Therefore, on \M relation $\nu=0$ holds and (\ref{bs-3}) becomes
\[R(X,Y)\xi=\kappa[\eta(Y)X-\eta(X)Y]+\mu[\eta(Y)AX-\eta(X)AY].\]
Thus, the structure vector field $\xi$ belongs to $(\kappa,\mu)$-nullity distribution, with $\kappa$,$\mu$ non-constant smooth functions.

\textbf{Case II:}  $\alpha^{2}+c=0$.

In this case the ambient space is $\mathbb{C}H^{n}$, $n\geq2$, and the above relation implies that $\alpha\neq0$. First suppose that $\lambda_{1}\neq\frac{\alpha}{2}$. Then relation (\ref{eq-B}) yields $\lambda_{2}=\frac{\alpha}{2}$. The inner product of relation (\ref{bs-3}) for $X=\varphi W$ and $Y=\xi$ with $W$ because of (\ref{eq-4}) and $\lambda_{2}=\frac{\alpha}{2}$ results in $\nu=0$. Thus, the structure vector field $\xi$ belongs to ($\kappa$,$\mu$)-nullity distribution, with $\kappa$,$\mu$ non-constant smooth functions. In previous section it has been proved that such real hypersurfaces in $M_{n}(c)$, $n\geq2$, do  not exist.

So $\lambda_{1}=\frac{\alpha}{2}$ is the only eigenvalue in $\mathbb{D}$ and \M is locally congruent to a horosphere. The inner product of relation (\ref{bs-3}) for $X=W$ and $Y=\xi$ with $\varphi W$ because of (\ref{eq-4}) and $\lambda_{1}=\frac{\alpha}{2}$ results in $\nu=0$. Thus, the structure vector field $\xi$ belongs to ($\kappa$,$\mu$)-nullity distribution, with $\kappa$,$\mu$ non-constant smooth functions. Therefore, because of Theorem \ref{tkn} the proof of Theorem \ref{tknm} is completed.

\begin{remark}
Let M be a real hypersurface in \MN,$n\geq2$, whose structure vector field $\xi$ belongs to ($\kappa,\mu,\nu$)-nullity distributions with $\kappa$,$\mu$,$\nu$ constants. Then following similar steps as in the proof of Theorem \ref{tknm} is proved

\begin{corollary}
Every real hypersurface in $M_{n}(c)$,$n\geq2$, whose structure vector $\xi$ belongs to ($\kappa,\mu,\nu$)-nullity distribution with $\kappa$,$\mu$,$\nu$ constants is Hopf.
\end{corollary}

Moreover, following similar steps to those of Hopf case in proof of Theorem \ref{tknm} it is concluded that

\begin{corollary}
Let M be a real hypersurface in $M_{n}(c)$,$n\geq2$, whose structure vector $\xi$ belongs to ($\kappa,\mu,\nu$)-nullity distribution with $\kappa$,$\mu$,$\nu$ constants. Then, $\xi$ belongs to ($\kappa,\mu)$-nullity distribution, with $\kappa$, $\mu$ constants and $\nu=0$.
\end{corollary}
\end{remark}
\section*{Open Problems}
\begin{itemize}
\item Firstly, it should be interesting to provide a complete classification of real hypersurfaces in non-flat complex space forms, whose structure vector field $\xi$ belongs to ($\kappa,\mu$)-nullity distribution with $\kappa$,$\mu$ constants.
\item Another interesting issue is to examine

\emph{if there exist real hypersurfaces in complex two-plane Grassmannians or complex hyperbolic two-plane Grassmannians (symmetric spaces of rank 2) whose structure vector field $\xi$ belongs to $\kappa$ or ($\kappa$,$\mu$) or ($\kappa$,$\mu$,$\nu$)-nullity distributions.}
\item Finally, it is known that the complex two-plane Grassmannians and the complex hyperbolic two-plane Grassmannians are equipped apart from the Kähler structure $J$ are also equipped with a quaternionic Kähler structure $\mathfrak{J}$ with local orthonormal basis $\{J_{1},J_{2},J_{3}\}$ which induces on \M an almost contact metric 3-structure ($\varphi_{i},\xi_{i},\eta_{i},g$), $i=1,2,3$ where $\xi_{i}=-J_{i}N$ and $N$ is the unit normal vector field on \M. Thus, another interesting issue is to examine

\emph{if there exist real hypersurfaces in complex two-plane Grassmannians or complex hyperbolic two-plane Grassmannians whose $\xi_{i}$, $i=1,2,3,$. belongs to $\kappa$- or ($\kappa$,$\mu$)- or ($\kappa$,$\mu$,$\nu$)-nullity distribution.}

Furthermore, on \M we can define the \textbf{($\kappa$,$\mu$,$\nu$)$_{i}$-nullity distribution} in the following way
\[N(\kappa, \mu, \nu):P\rightarrow N_{P}(\kappa, \mu, \nu),\;\;\mbox{where $\kappa$, $\mu$, $\nu$ are functions}\]
 and  $N_{P}(\kappa, \mu, \nu)$ is given by
\[N_{P}(\kappa, \mu ,\nu)=\{Z\;\in\;T_{P}M:R(X,Y)Z=\kappa[\eta_{i}(Y)X-\eta_{i}(X)Y]+\mu[\eta_{i}(Y)AX-\eta_{i}(X)AY]\]
\[+\nu[\eta_{i}(Y)\varphi_{i} AX-\eta_{i}(X)\varphi_{i} AY]\},\;\;\mbox{for any $X$, $Y$ $\in$ $T_{P}M$ and $i=1,2,3$}.\]
So the following questions raises naturally

\emph{Are there real hypersurfaces in the above spaces whose structure vector field $\xi$ belongs to  ($\kappa$,$\mu$,$\nu$)$_{i}$-nullity distribution?}

\emph{Are there real hypersurfaces in the above spaces whose $\xi_{i}$ belongs to  ($\kappa$,$\mu$,$\nu$)$_{i}$-nullity distribution?}
\end{itemize}

\section*{Acknowledgements}
The author would like to thank Prof. Ph. J. Xenos for his comments on the paper.

\scriptsize{\textsc{\hspace{-15pt}K. Panagiotidou, Faculty of Mathematics and Engineering Sciences, Hellenic Military Academy, Vari, Attiki, Greece}}\\
\textsc{e-mail}: konpanagiotidou@gmail.com

 \end{document}